\documentclass[12pt,a4paper,reqno]{amsart}
\usepackage{amssymb}
\usepackage{amscd}
\usepackage{enumerate}
\usepackage{graphicx}
\usepackage{siunitx}
\usepackage{tikz-cd}
\usepackage{bm}
\numberwithin{equation}{section}

\usepackage{mathtools}
\usepackage[tableposition=top]{caption}
\usepackage{booktabs,dcolumn}

\DeclareFontFamily{OT1}{rsfs}{}
\DeclareFontShape{OT1}{rsfs}{n}{it}{<-> rsfs10}{}
\DeclareMathAlphabet{\mathscr}{OT1}{rsfs}{n}{it}

\theoremstyle{plain}

\newtheorem{theorem}{Theorem}[section]
\newtheorem{proposition}[theorem]{Proposition}

\theoremstyle{definition}

\newtheorem{remark}[theorem]{Remark}

\newcommand\R{\mathbb{R}}

\begin{document}

\title[Perfectly packing a square by squares]{Perfectly packing a square by squares of nearly harmonic sidelength}

\author{Terence Tao}
\address{UCLA Department of Mathematics, Los Angeles, CA 90095-1555.}
\email{tao@math.ucla.edu}



\begin{abstract}  A well known open problem of Meir and Moser asks if the squares of sidelength $1/n$ for $n \geq 2$ can be packed perfectly into a square of area $\sum_{n=2}^\infty \frac{1}{n^2} = \frac{\pi^2}{6}-1$.    In this paper we show that for any $1/2 < t < 1$, and any $n_0$ that is sufficiently large depending on $t$, the squares of sidelength $n^{-t}$ for $n \geq n_0$ can be packed perfectly into a square of area $\sum_{n=n_0}^\infty \frac{1}{n^{2t}}$.  This was previously known (if one packs a rectangle instead of a square) for $1/2 < t \leq 2/3$ (in which case one can take $n_0=1$).
\end{abstract}

\maketitle


\section{Introduction}

A \emph{packing} by rectangles\footnote{In this paper all rectangles and squares are understood to have sides parallel to the coordinate axes, and to be topologically closed.} of a region $\Omega \subset \R^2$ is a finite or countably infinite family of rectangles in $\Omega$ with disjoint interiors.  We say that the packing is \emph{perfect} if the rectangles cover $\Omega$ up to null sets.  Note that this forces the Lebesgue measure $m(\Omega)$ of $\Omega$ to equal the sum $\sum_{n=1}^\infty \mathrm{area}(R_n)$ of the areas of the rectangles.

Meir and Moser \cite{meir} posed the question of whether rectangles of dimensions $\frac{1}{n} \times \frac{1}{n+1}$ for $n \geq 1$ can perfectly pack the unit square $[0,1]^2$, as well as the very similar question of whether squares of sidelength $\frac{1}{n}$ for $n \geq 2$ can perfectly pack a square of area $\sum_{n=2}^\infty \frac{1}{n^2} = \frac{\pi^2}{6}-1$.  These questions remain open; see for instance \cite{croft}, \cite[Chapter 3]{brass} for further discussion.  As one measure of partial progress towards these results, Paulhus \cite{paulhus} showed\footnote{As pointed out in \cite{joos-0}, some of the lemmas in this paper were not proven correctly, but the gaps in this paper were recently repaired in \cite{gj}.} that one could pack rectangles of dimensions $\frac{1}{n} \times \frac{1}{n+1}$ for $n \geq 1$ into a square of area $1 + \frac{1}{10^9+1}$, and squares of sidelength $1/n$ for $n \geq 2$ into a rectangle of area $\frac{\pi^2}{6} - 1 + \frac{1}{1244918662}$.

Another direction in which partial progress has been made is to consider whether, for any $t > 1/2$, squares of sidelength $n^{-t}$ for $n \geq 1$ can perfectly pack a square or rectangle of area $\sum_{n=1}^\infty \frac{1}{n^{2t}}$ (which is finite when $t>1/2$).  The goal is then to get $t$ as close as possible to $1$, to address the second question of Meir and Moser posed above.  Recently an affirmative answer to this question was given in the range $1/2 < t \leq 2/3$ by Januszewski and Zielonka \cite{jz}, building upon previous work in \cite{chalcraft}, \cite{wastlund}, \cite{joos}, as well as a packing algorithm in the previously mentioned paper \cite{paulhus}.  In this note we extend the range of $t$ to almost reach the value $t=1$ corresponding to the question of Meir and Moser, at the expense of removing the first few squares in the sequence:

\begin{theorem}\label{main}  Let $1/2 < t < 1$, and suppose that $n_0$ is a natural number that is sufficiently large depending on $t$.  Then squares of sidelength $n^{-t}$ for $n \geq n_0$ can perfectly pack a square of area $\sum_{n=n_0}^\infty \frac{1}{n^{2t}}$.
\end{theorem}

As a corollary, for every $1/2 < t < 1$, the squares of sidelength $n^{-t}$ for $n \geq 1$ can perfectly pack a finite union of squares; this latter claim was first established in the range $1/2 < t \leq 5/6$ in \cite{chalcraft}, and extended to the range $1/2 < t < 2/3$ in \cite{wastlund}.

The strategy of proof is similar to that in the previous works \cite{chalcraft}, \cite{wastlund}, \cite{joos}, \cite{jz}, in which one performs a recursive algorithm to pack the first few squares $n^{-t}$, $n_0 \leq n < n_1$ into a square of the indicated area, with the remaining space being described by a union of a family ${\mathcal R}_{n_1}$ of rectangles which have a certain controlled size.  In previous algorithms, the total perimeter of this family ${\mathcal R}_{n_1}$ was comparable to the total perimeter $\sum_{n=n_0}^{n_1-1} \frac{4}{n^t}$ of the squares that one had already packed, and thus (for large $n_1$) also comparable to $n_1^t$ times the total area $\sum_{n=n_1}^\infty \frac{1}{n^{2t}}$ of the remaining rectangles.
 It is this relationship between the total perimeter and total area of ${\mathcal R}_{n_1}$ that prevents $t$ from getting too close to $1$, as otherwise one could not eliminate the possibility that all remaining rectangles in ${\mathcal R}$ had width less than $n_1^{-t}$, thus preventing one from continuing the packing.  By arranging the squares in near-lattice formations, we are able (for $n_0$ large enough) to make the total perimeter of ${\mathcal R}_{n_1}$ significantly smaller than the perimeter of the squares that one has already packed, and thus significantly smaller than $n_1^t$ times the total area of ${\mathcal R}_{n_1}$; this will allow us to take $t$ arbitrarily close to $1$.  Unfortunately the argument does not seem to extend to the critical case $t=1$ (or to the supercritical cases $t>1$).

We remark that the same argument (with minor notational changes) would also allow one to pack rectangles of dimensions $n^{-t} \times (n+1)^{-t}$ for $n \geq n_0$ perfectly into a square of area $\sum_{n=n_0}^\infty \frac{1}{n^t(n+1)^t}$; we leave the details of this modification to the interested reader.  The quantity $n_0$ could be calculated explicitly as a function of $t$, but we have not attempted to optimize this quantity.  In principle, one could combine the arguments here with some initial packing of the first $n_0$ squares, located for instance by computer search, in order to be able to replace $n_0$ by $1$ for certain values of $t$ that are sufficiently far from $1$, but we will not attempt to do so here.

\subsection{Acknowledgments}

The author is supported by NSF grant DMS-1764034 and by a Simons Investigator Award.  We thank Rachel Greenfeld, Jose Madrid, Keiju Sono and anonymous commentators on the author's blog for corrections.

\subsection{Data sharing statement}

No datasets were generated in this paper.

\section{Initial reductions}

Throughout this paper we fix the parameter $1/2 < t < 1$, and then introduce the exponent
$$ \delta \coloneqq 1-t$$
Note that because we are in the regime $1/2 < t < 1$, we have $0 < \delta < 1$ and
\begin{equation}\label{delta-small}
t + \delta t < 1.
\end{equation}
In fact, these are the only two properties of $\delta$ that we will need in the sequel.  We will use this exponent $\delta$ to define a certain technical modification of the concept of the total perimeter of a family of rectangles.

We adopt the asymptotic notation $X = O(Y)$, $X \ll Y$, or $Y \gg X$ to denote the estimate $|X| \leq C_{t} Y$ for some constant $C_{t}$ that is allowed to depend only on $t$ (or equivalently, on $\delta$);  in particular, these constants will be independent of the parameters $M$ or $N_0$ that we shall shortly introduce.  We write $X \asymp Y$ for $X \ll Y \ll X$.

Next, we select two large parameters:

\begin{itemize}
\item We pick a natural number $M$ which is sufficiently large depending on $\delta,t$.  (One can for instance take $M \coloneqq \lfloor \exp(C/\delta) \rfloor$ for a suitably large absolute constant $C$.)  Roughly speaking, we will pack our squares in groups of cardinality $\asymp M^2$ at a time, arranged into approximate lattices with $\asymp M$ squares in each row and column.
\item Finally, we pick a number $N_0$ that is sufficiently large depending on $M, \delta, t$.  (For instance, one can check that $N_0 \coloneqq M^{10/\delta}$ would work in the arguments below, though this choice is far from best possible.)  This will be our lower bound for the parameter $n_0$ in Theorem \ref{main}; in particular, $n_0$ will be far larger than $M$ or $M^2$.
\end{itemize}

Given a rectangle $R$, we define the \emph{width} $w(R)$ to be the smaller of the two sidelengths, and the \emph{height} $h(R)$ to be the larger of the two sidelengths (with $w(R)=h(R)$ when $R$ is a square), thus the area $m(R)$ is equal to $w(R) h(R)$.  Given a finite family ${\mathcal R}$ of rectangles with disjoint interiors, we can thus define the \emph{total area}
$$ \mathrm{area}({\mathcal R}) \coloneqq \sum_{R \in {\mathcal R}} w(R) h(R)$$
and
\emph{unweighted total perimeter}
$$ \mathrm{perim}({\mathcal R}) \coloneqq \sum_{R \in {\mathcal R}} 2(w(R)+h(R)) \asymp \sum_{R \in {\mathcal R}} h(R).$$
For technical reasons we will often work instead with the \emph{weighted total perimeter} 
$$ \mathrm{perim}_\delta({\mathcal R}) \coloneqq \sum_{R \in {\mathcal R}} w(R)^\delta h(R).$$
One should think of this weighted total perimeter as a slight modification of the unweighted total perimeter, in which narrower rectangles are given slightly less weight than wider rectangles.  This modification is convenient for technical induction purposes; our algorithms will at one point replace a wide rectangle with several narrower rectangles, with a favorable control on the weighted total perimeter of the latter, despite having unfavorable control on the unweighted total perimeter.

In previous literature, proofs of results such as Theorem \ref{main} were given by detailing a specific recursive algorithm for generating the desired packing, and then verifying that the algorithm produced a packing with all the required properties.  Here we will arrange the argument slightly differently\footnote{See however Remark \ref{algo} below.} by using induction instead of recursion, and more precisely by using a downward induction to establish the following more technical proposition, that allows us to perfectly pack any family of rectangles that has well controlled weighted total perimeter (and also obeys some other minor conditions), and which easily implies Theorem \ref{main}:

\begin{proposition}[Perfectly packing some families of rectangles]\label{nn1}  Let $n_{\max} \geq n_0 \geq N_0$, and suppose that ${\mathcal R}$ is a finite family of rectangles with disjoint interiors, with total area 
\begin{equation}\label{total-area}
\mathrm{area}({\mathcal R}) = \sum_{n=n_0}^\infty \frac{1}{n^{2t}},
\end{equation}
 and obeying the weighted total perimeter bound
\begin{equation}\label{wtr}
 \mathrm{perim}_\delta({\mathcal R})  \leq M^{-1 + \frac{\delta}{2}} \sum_{n=1}^{n_0-1} \frac{1}{n^{t+\delta t}}
\end{equation}
and the crude height bound
\begin{equation}\label{hr}
 \sup_{R \in {\mathcal R}} h(R) \leq 1.
\end{equation}
Then one can pack $\bigcup_{R \in {\mathcal R}} R$ by squares of sidelength $n^{-t}$ for $n_0 \leq n < n_{\max}$.
\end{proposition}

Indeed, if $n_0 \geq N_0$ and we take ${\mathcal R}$ to consist solely of a square $S$ of area $\sum_{n=n_0}^\infty \frac{1}{n^{2t}}$, then $S$ has sidelength $O( n_0^{1/2-t} )$ (here we use the hypothesis $t>1/2$), and hence
$$ \mathrm{perim}_\delta({\mathcal R}) \ll n_0^{(1/2-t)(1+\delta)}.$$
On the other hand, from \eqref{delta-small} we have
$$ \sum_{n=1}^{n_0-1} \frac{1}{n^{t+\delta t}} \gg n_0^{1-t-\delta t} = n_0^{\frac{1-\delta}{2}} n_0^{(1/2-t)(1+\delta)}.$$
Since $n_0 \geq N_0$ and $N_0$ is sufficiently large depending on $M,\delta,t$, we conclude that the condition \eqref{wtr} holds.  Also it is clear that $S$ has height at most $1$.  Applying Proposition \ref{nn1}, we conclude that we can pack $S$ by the squares of sidelength $n^{-t}$ for $n_0 \leq n < n_{\max}$ for any $n_{\max}$.  Sending $n_{\max} \to \infty$ and using a standard compactness argument (see e.g., \cite{martin}) we can then pack $S$ by squares of sidelength $n^{-t}$ for $n \geq n_0$, which is then a perfect packing by comparison of areas.  Theorem \ref{main} follows.

The key step in establishing Proposition \ref{nn1} will be to prove the following assertion.

\begin{proposition}[Efficiently packing a small rectangle of bounded eccentricity]\label{nn-eff}  Let $n_0 \geq N_0$, and suppose that $R$ is a rectangle whose dimensions $w(R), h(R)$ obey the inequalities
\begin{equation}\label{Ma}
 M n_0^{-t} \leq w(R) \leq h(R) \leq 3 M n_0^{-t}.
\end{equation}
Then one can find $n'_0 \geq n_0$ with $n'_0 - n_0 \asymp M^2$ and a perfect packing of $R$ by the squares of sidelength $n^{-t}$ for $n_0 \leq n < n'_0$, together with an additional finite family ${\mathcal R}$ of rectangles with disjoint interiors and widths $O(n_0^{-t})$, obeying the unweighted total perimeter bound
\begin{equation}\label{hrr}
\mathrm{perim}({\mathcal R}) \ll M n_0^{-t}.
\end{equation}
\end{proposition}

The point here is that the unweighted total perimeter of the rectangles ${\mathcal R}$ is only $O(M n_0^{-t})$, as compared against the unweighted total perimeter of the squares of sidelength $n^{-t}$ for $n_0 \leq n < n'_0$ which is comparable to $M^2 n_0^{-t}$.  This gain of $O(M^{-1})$ is superior to the factor of $M^{-1+\delta/2}$ which appears in \eqref{wtr}, which in turn is superior to the factor $M^{-1+\delta}$ which is what would be needed to ensure the condition \eqref{Ma} is satisfied for certain rectangles $R_i$ that we will construct shortly.

We prove Proposition \ref{nn-eff} in the next section.  Assuming it for now, we conclude the proof of Proposition \ref{nn1} and hence Theorem \ref{main}.  We fix $n_{\max}$ and perform a downward induction on $n_0$.  Proposition \ref{nn1} is trivially true for $n_0=n_{\max}$, so suppose that $n_0<n_{\max}$ and that the claim has already been proven for larger values of $n_0$.  From \eqref{wtr} and \eqref{delta-small} we have
$$
\sum_{R \in {\mathcal R}} w(R)^\delta h(R) \ll M^{1 + \frac{\delta}{2}} n_0^{1-t-\delta t}.$$
On the other hand, from \eqref{total-area} we have
$$
\sum_{R \in {\mathcal R}} w(R) h(R) \gg n_0^{1-2t}.$$
From the pigeonhole principle, we conclude that there exists $R \in {\mathcal R}$ with
$$ w(R)^{1-\delta} \gg \frac{n_0^{1-2t}}{M^{-1 + \frac{\delta}{2}} n_0^{1-t-\delta t}}$$
which simplifies (using $(1-\delta/2) / (1-\delta) > (1+\delta/2)$) to
$$ w(R) \gg M^{1+\frac{\delta}{2}} n_0^{-t}.$$
Of course this implies
\begin{equation}\label{har}
 h(R) \geq w(R) \geq 2M n_0^{-t}.
\end{equation}
We can then partition $R$ into a rectangle $R_0$ of dimensions $(w(R)-Mn_0^{-t}) \times h(R)$ and a rectangle $R_*$ of dimensions $Mn_0^{-t} \times h(R)$.  By cutting off squares of sidelength $Mn_0^{-t}$ from $R_*$ until the height of the remaining rectangle dips below $2Mn_0^{-t}$, we see from \eqref{har} that one can partition $R_*$ into rectangles $R_1,\dots,R_m$ of dimensions $Mn_0^{-t} \times h(R_i)$ with
$$ Mn_0^{-t} \leq h(R_i) < 2Mn_0^{-t}$$
for $i=1,\dots,m$, and
$$ \sum_{i=1}^m h(R_i) = h(R).$$
From \eqref{hr} we conclude in particular the crude upper bound
\begin{equation}\label{crude-upper}
m \leq n_0^{t}
\end{equation}
and we have the perfect packing
\begin{equation}\label{rro}
 R = R_0 \cup R_* = R_0 \cup R_1 \cup \dots \cup R_m.
\end{equation}

Applying Proposition \ref{nn-eff} $m$ times, we can then find natural numbers
$$ n_0 = n'_0 < n'_1 < \dots < n'_m$$
with 
\begin{equation}\label{njump}
n'_{i+1} - n'_i \asymp M^2
\end{equation}
for all $0 \leq i \leq m-1$, which by \eqref{crude-upper} implies in particular that 
\begin{equation}\label{ni}
n_0 \leq n'_i \leq 1.001 n_0
\end{equation}
(say) for all $0 \leq i \leq m$, and a perfect packing of each $R_i$, $i=1,\dots,m$ by squares of sidelength $n^{-t}$ for $n'_{i-1} \leq n < n'_i$, together with an additional family ${\mathcal R}_i$ of rectangles of disjoint interiors, widths $O( n_0^{-t})$, and with
\begin{equation}\label{harm}
\mathrm{perim}({\mathcal R}_i) \ll M n_0^{-t}.
\end{equation}

If we then define the new family of rectangles
$$ {\mathcal R}' \coloneqq ({\mathcal R} \backslash \{R\}) \cup \{R_0\} \cup \bigcup_{i=1}^m {\mathcal R}_i$$
then we see that the rectangles in ${\mathcal R}'$ have disjoint interiors, and $\bigcup_{R' \in {\mathcal R}} R'$ is perfectly packed by squares of sidelength $n^{-t}$ for $n_0 \leq n < n'_m$, together with the rectangles in ${\mathcal R}'$.  If $n'_m \geq n_{\max}$ then we are now done, so assume that $n'_m < n_{\max}$.  We compute (using $w(R_0) \leq w(R)$, $h(R_0) = h(R)$, \eqref{harm}, \eqref{ni}, \eqref{njump}, and \eqref{wtr} in turn)
\begin{align*}
\mathrm{perim}_\delta({\mathcal R}') &= \mathrm{perim}_\delta({\mathcal R}) - w(R)^\delta h(R) + w(R_0)^\delta h(R_0) + \sum_{i=1}^m \sum_{R' \in {\mathcal R}_i} w(R')^\delta h(R') \\
&\leq \mathrm{perim}_\delta({\mathcal R}) + \sum_{i=1}^m O( n_0^{-\delta t} \mathrm{perim}({\mathcal R}_i) ) \\
&\leq \mathrm{perim}_\delta({\mathcal R}) + M \sum_{i=1}^m O( n_0^{-t-\delta t} ) \\
&= \mathrm{perim}_\delta({\mathcal R}) + M^{-1} \sum_{i=1}^m O\left( \sum_{n=n'_i}^{n'_{i+1}-1} \frac{1}{n^{t+\delta t}} \right) \\
&= \mathrm{perim}_\delta({\mathcal R}) + M^{-1} O\left( \sum_{n=n_0}^{n'_m-1} \frac{1}{n^{t+\delta t}} \right) \\
&\leq M^{-1 + \frac{\delta}{2}} \sum_{n=1}^{n'_m-1} \frac{1}{n^{t+\delta t}};
\end{align*}
that is to say, ${\mathcal R}'$ obeys the condition \eqref{wtr} (with $n_0$ replaced by $n'_m$).  Also, the total area of ${\mathcal R}'$ can be computed to be
$$ \mathrm{area}({\mathcal R}') = \mathrm{area}({\mathcal R}) - \sum_{n = n_0}^{n'_m-1} \frac{1}{n^{2t}} = \sum_{n=n'_m}^\infty \frac{1}{n^{2t}}$$
and from \eqref{hr} we easily see that all rectangles in ${\mathcal R}'$ have height at most $1$.  Thus by induction hypothesis, we can pack $\bigcup_{R' \in {\mathcal R}'} R'$ by squares of sidelength $n^{-t}$ for $n'_m \leq n < n_{\max}$.  This gives the desired packing of ${\mathcal R}$ by squares of sidelength $n^{-t}$ for $n_0 \leq n < n_{\max}$, closing the induction.

It remains to establish Proposition \ref{nn-eff}.  This is the purpose of the next section.

\begin{remark}\label{algo} The above analysis can be converted into the following algorithm for constructing the perfect packing in Theorem \ref{main}:
\begin{enumerate}
\item Select a sufficiently large natural number $M$, initialize $n_0$ to be the quantity in Theorem \ref{main}, and let ${\mathcal R}$ consist of a single square $S$ of sidelength $\sum_{n=n_0}^\infty \frac{1}{n^{2t}}$.
\item Let $R$ be a rectangle in ${\mathcal R}$ of maximal width $w(R)$, and perform the subdivision \eqref{rro} of $R$ into rectangles $R_0,R_1,\dots,R_m$ as indicated above. (This assumes that $w(R) \geq 2Mn_0^{-t}$; if this is not the case, terminate with error.)
\item For each $i=1,\dots,m$ in turn, apply Proposition \ref{nn-eff} to $R_i$ to subdivide that rectangle into squares of sidelength $n^{-t}$ for $n_0 \leq n < n'_0$, together with an additional family of rectangles ${\mathcal R}_i$; then replace $n_0$ with $n'_0$ and continue iterating in $i$.
\item Replace the rectangle $R$ in ${\mathcal R}$ by $R_0$ together with the rectangles in ${\mathcal R}_1 \cup \dots \cup {\mathcal R}_m$, then return to Step 2.
\end{enumerate}
The above analysis then ensures (for $n_0$ large enough) that this algorithm never terminates and produces a perfect packing of the original square $S$.
\end{remark}

\section{Efficiently packing a small rectangle of bounded eccentricity}

We now prove Proposition \ref{nn-eff}. Without loss of generality we may take $R$ to be the rectangle
$$R = [0, w(R)] \times [0, h(R)].$$
From \eqref{Ma} we may find natural numbers
$$ M \leq M_1 \leq M_2 < 3M$$
such that
$$ M_1 n_0^{-t} \leq w(R) < (M_1+1) n_0^{-t}$$
and
$$ M_2 n_0^{-t} \leq h(R) < (M_2+1) n_0^{-t}.$$
We will now take $n'_0 \coloneqq n_0 + M_1 M_2$, then clearly $n'_0 - n_0 \asymp M^2$.  We index the set $\{ n: n_0 \leq n < n_0 + M_1 M_2\}$ ``lexicographically'' as $\{ n_{i,j}: 0 \leq i < M_1; 0 \leq j < M_2\}$, where
$$ n_{i,j} \coloneqq n_0 + j M_1 + i.$$
Our task is then to perfectly pack $R$ by $M_1 M_2$ squares $S_{i,j}$ of sidelength $1/n_{i,j}^t$ for $0 \leq i < M_1$ and $0 \leq j < M_2$, together with some additional finite family ${\mathcal R}$ of rectangles with disjoint interiors and heights $O(n_0^{-t})$ obeying \eqref{hrr}.

To motivate the construction, suppose temporarily that the squares $S_{i,j}$ were required to have sidelength $1/n_0^t$ instead of $1/n_{i,j}^t$. Then we could simply use the lattice packing
\begin{equation}\label{lat}
 S_{i,j} \coloneqq [i n_0^{-t}, (i+1)n_0^{-t}] \times [j n_0^{-t}, (j+1)n_0^{-t}]
\end{equation}
for $0 \leq i < M_1, 0 \leq j < M_2$, as these squares perfectly pack the rectangle
$$ [0, M_1 n_0^{-t}] \times [0, M_2 n_0^{-t}]$$
and the remaining portion of the original rectangle $R$ can then be perfectly packed by the two rectangles
$$ [0, M_1 n_0^{-t}] \times [M_2 n_0^{-t}, h(R)]$$
and
$$ [M_1 n_0^{-t}, w(R)] \times [0, h(R)]$$
which have widths $O(n_0^{-t})$ and heights $O(M n_0^{-t})$ (and thus perimeters $O(M n_0^{-t})$), giving the claim.

\begin{figure} [t]
\centering
\includegraphics{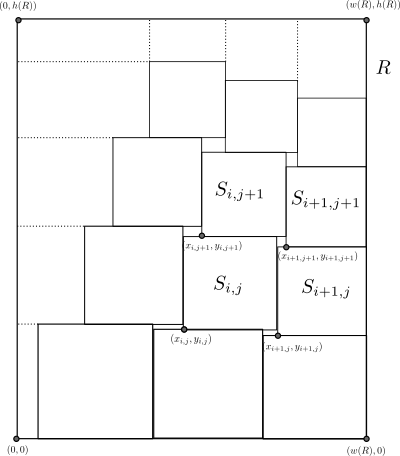}
\caption{A rectangle $R$ (with $M_1=3$ and $M_2=4$), which is perfectly packed by $M_1 M_2 = 12$ squares $S_{i,j}$ with $0 \leq i < 3$ and $0 \leq j < 4$ (the square $S_{i,j}$ depicted is for $(i,j)=(1,1)$), together with $(M_1-1)(M_2-1)=6$ small rectangles of the form \eqref{surround} between the squares $S_{i,j}$, $M_2=4$ rectangles of the form \eqref{m2} on the left side of $R$, $M_1=3$ rectangles of the form \eqref{m1} on the upper side of $R$, and one rectangle \eqref{m0} on the upper left of $R$.  This becomes a reasonably efficient packing of the rectangle $R$ by squares once $M$ (and hence $M_1,M_2$) gets large, and $n_0$ is extremely large compared to $M$.}
\label{fig:squish}
\end{figure}

In our actual problem, the squares $S_{i,j}$ are slightly smaller, being required to have sidelength $n^{-t}$ instead of $1/n_0^t$.  If one attempts to position the bottom left corners of the $S_{i,j}$ in the same location $(i n_0^{-t}, j n_0^{-t})$ as in the lattice packing \eqref{lat}, thus
$$ S_{i,j} \coloneqq [i n_0^{-t}, in_0^{-t} + n_{i,j}^{-t}] \times [j n_0^{-t}, jn_0^{-t} + n_{i,j}^{-t}]$$
then this would still form a packing of the rectangle $R$, but there would now be a large number of gaps between the squares, necessitating ${\mathcal R}$ to consist of something like $\asymp M^2$ rectangles of perimeter $\asymp n_0^{-t}$ each, which would not give the desired bound \eqref{hrr}.  However, it is possible to close most of these gaps by sliding the squares $S_{i,j}$ closer together, thus reducing the perimeter of ${\mathcal R}$ substantially.  More precisely, our actual construction of the $S_{i,j}$ will take the form
$$ S_{i,j} \coloneqq [x_{i,j}, x_{i,j} + n_{i,j}^{-t}] \times [y_{i,j}, y_{i,j} + n_{i,j}^{-t}]$$
where
$$ x_{i,j} \coloneqq w(R) - \sum_{i'=i}^{M_1-1} n_{i',j}^{-t}$$
and
$$ y_{i,j} \coloneqq \sum_{j'=0}^{j-1} n_{i,j'}^{-t};$$
see Figure \ref{fig:squish}.
Note from the mean value theorem, the triangle inequality, and the hypothesis $n_0 \geq N_0$ that
\begin{equation}\label{x-approx}
 x_{i,j} = w(R) - (M_1-i+1 + O( M^3/N_0) ) n_0^{-t} 
\end{equation}
and
\begin{equation}\label{y-approx}
 y_{i,j} = (j + O( M^3/N_0) ) n_0^{-t}.
\end{equation}
Thus, up to errors of $O(\frac{M^3}{N_0} n_0^{-t})$, the points $(x_{i,j}, y_{i,j})$ are arranged in a lattice of spacing $n_0^{-t}$.  Note that for any $0 \leq i < M_1$ and $0 \leq j < M_2$ we have
$$ 0 \leq w(R) - M_1 n_0^{-t} \leq x_{i,j} \leq x_{i,j} + n_{i,j}^{-t} \leq w(R)$$
and
$$ 0 \leq y_{i,j} \leq y_{i,j} + n_{i,j}^{-t} \leq M_2 n_0^{-t} \leq h(R)$$
and so all the squares $S_{i,j}$ are contained in $R$.  Next, for any $0 \leq i, i' < M_1$ and $0 \leq j,j' < M_2$ with $(i,j) \neq (i',j')$, we argue that the squares $S_{i,j}, S_{i',j'}$ have disjoint interiors as follows.

\begin{itemize}
\item If $i' < i$ and $j' \geq j$, then $y_{i',j'} + n_{i',j'}^{-t} \leq y_{i,j}$, and hence the interior of $S_{i',j'}$ lies below the interior of $S_{i,j}$, giving disjointness.  Similarly if $i < i'$ and $j \geq j'$.
\item If $j' < j$ and $i' \leq i$, then $x_{i',j'} + n_{i',j'}^{-t} \leq x_{i,j}$, and hence the interior of $S_{i',j'}$ lies to the left of the interior of $S_{i,j}$, giving disjointness.  Similarly if $j < j'$ and $i \leq i'$.  This covers all the possible cases for $i,j,i',j'$.
\end{itemize}

If $0 \leq i < M_1-1$ and $0 \leq j < M_2-1$, then (using \eqref{x-approx}, \eqref{y-approx} as necessary) we have the relations
\begin{align*}
x_{i+1,j} &= x_{i,j} + n_{i,j}^{-t} \\
y_{i,j} &< y_{i+1,j} + n_{i+1,j}^{-t} < y_{i,j} + n_{i,j}^{-t} \\
x_{i,j} &< x_{i,j+1} < x_{i,j} + n_{i,j}^{-t} \\
y_{i,j+1} &= y_{i,j} + n_{i,j}^{-t}\\
x_{i+1,j+1} &= x_{i,j+1} + n_{i,j+1}^{-t} \\
y_{i+1,j+1} &< y_{i,j+1} < y_{i+1,j+1} + n_{i+1,j+1}^{-t} \\
x_{i+1,j} &< x_{i+1,j+1} < x_{i+1,j} + n_{i+1,j}^{-t} \\
y_{i+1,j+1} &= y_{i+1,j} + n_{i+1,j}^{-t}
\end{align*}
(see Figure \ref{fig:squish}).  As a consequence, the squares $S_{i,j}, S_{i+1,j}, S_{i,j+1}, S_{i+1,j+1}$ surround the rectangle
\begin{equation}\label{surround}
[x_{i,j+1}, x_{i+1,j+1}] \times [y_{i+1,j+1}, y_{i,j+1}]
\end{equation}
which by \eqref{x-approx}, \eqref{y-approx} has width and height $O( \frac{M^3}{N_0} n_0^{-t} )$, and hence perimeter $O( \frac{M^3}{N_0} n_0^{-t} )$ also.

From Figure \ref{fig:squish} we now see that the rectangle $R$ can be packed by the squares $S_{i,j}$ for $0 \leq i < M_1, 0 \leq j < M_2$ together with the rectangles \eqref{surround} for $0 \leq i < M_1-1, 0 \leq j < M_2-1$, as well as the additional rectangles
\begin{equation}\label{m2}
[0, x_{0,j}] \times [y_{0,j}, y_{0,j} + n_{0,j}^{-t}]
\end{equation}
for $0 \leq j < M_2$, the rectangles
\begin{equation}\label{m1}
[x_{i,M_2-1}, x_{i,M_2-1} + n_{i,M_2-1}^{-t}] \times [y_{i,M_2-1}+n_{i,M_2-1}^{-t}, h(R)]
\end{equation}
for $0 \leq i < M_1$, and the rectangle
\begin{equation}\label{m0}
 [0, x_{0,M_2-1}] \times [y_{0,M_2-1} + n_{0,M_2-1}^{-t}, h(R)].
\end{equation}
All of these rectangles have width and height $O(n_0^{-t})$, thanks to \eqref{x-approx}, \eqref{y-approx}, and hence perimeter $O(n_0^{-t})$ also.  Collecting these rectangles into a family ${\mathcal R}'$, we see that
$$
\mathrm{perim}({\mathcal R}) \ll M^2 \times \frac{M^3}{N_0} n_0^{-t} + M \times n_0^{-t}
$$
which gives \eqref{hrr} since $N_0$ is large compared with $M$.  The claim follows.


\begin{thebibliography}{10}

\bibitem{brass}
P. Brass, W. Moser, J. Pach, Research problems in discrete geometry. Springer, New York, 2005.

\bibitem{chalcraft}
A. Chalcraft, \emph{Perfect square packings}, J. Combin. Theory Ser. A \textbf{92} (2000), 158--172.

\bibitem{croft}
H. T. Croft, K. J. Falconer, R. K. Guy, \emph{Unsolved Problems in Geometry}, pp. 112–113, Springer-Verlag, New York, 1991.

\bibitem{gj}
P. Grzegorek, J. Januszewski, \emph{A note on three Moser's problems and two Paulhus' lemmas}, J. Combin. Theory Ser. A \textbf{162} (2019), 222--230.

\bibitem{jz}
J. Januszewski, {\L}. Zielonka, \emph{A note on perfect packing of squares and cubes}, Acta Math. Hungar. 163 (2021), no. 2, 530--537.

\bibitem{joos-0}
A. Jo\'os, \emph{On packing of rectangles in a rectangle}, Discrete Math. \textbf{341} (2018), no. 9, 2544--2552.

\bibitem{joos}
A. Jo\'os, \emph{Perfect packing of cubes}, Acta Math. Hungar., \textbf{156} (2018), 375--384.

\bibitem{joos-2}
A. Jo\'os, \emph{Perfect square packings}, preprint.

\bibitem{martin}
G. Martin, \emph{Compactness theorems for geometric packings}, J. Combin. Theory Ser. A \textbf{97} (2002), no. 2, 225--238.

\bibitem{meir}
A. Meir, L. Moser, \emph{On packing of squares and cubes}, J. Combin. Theory \textbf{5} (1968), 126--134.

\bibitem{paulhus}
M. M. Paulhus, \emph{An algorithm for packing squares}, J. Combin. Theory Ser. A \textbf{82} (1998), 147--157.

\bibitem{wastlund}
J. W\"astlund, \emph{Perfect packings of squares using the stack-pack strategy}, Discrete Comput. Geom. 29 (2003), no. 4, 625--631.

\end{thebibliography}
\end{document}